\documentclass{siamart251216}

\usepackage{amssymb,mathtools}
\usepackage{bm}
\usepackage{enumitem}
\usepackage{algpseudocode}
\usepackage{microtype}
\usepackage{needspace}
\usepackage{tikz}
\usepackage{stmaryrd}
\usetikzlibrary{arrows.meta,positioning}

\newsiamremark{remark}{Remark}
\newsiamremark{example}{Example}
\newsiamthm{assumption}{Assumption}
\crefname{assumption}{assumption}{assumptions}
\Crefname{assumption}{Assumption}{Assumptions}
\crefname{remark}{remark}{remarks}
\Crefname{remark}{Remark}{Remarks}
\crefname{example}{example}{examples}
\Crefname{example}{Example}{Examples}
\AddToHook{env/assumption/begin}{\crefalias{theorem}{assumption}}
\AddToHook{env/remark/begin}{\crefalias{theorem}{remark}}
\AddToHook{env/example/begin}{\crefalias{theorem}{example}}

\newcommand{\A}{\mathbb A}
\newcommand{\dd}{\mathrm{d}}
\newcommand{\osc}{\operatorname{osc}}
\newcommand{\card}{\#}
\newcommand{\en}[1]{\lVert\nabla(#1)\rVert_{L^2(\Omega)}}

\headers{Rate optimality of AFEM with maximum marking}{Y. Li}
\title{Rate Optimality of AFEM with Maximum Marking}
\author{Yuwen Li\thanks{School of Mathematical Sciences, Zhejiang University,
866 Yuhangtang Road, Hangzhou 310058, Zhejiang, People's Republic of China
(\email{liyuwen@zju.edu.cn}).}}

\begin{document}
\maketitle

\begin{abstract}
For Poisson's equation, we prove optimal algebraic convergence rate for an adaptive finite element method (AFEM) driven by standard maximum marking in arbitrary space dimension, for every marking parameter in $\lparen 0,1\rbrack$. Previously, instance optimality of a modified maximum marking was proved in {\sc L.~Diening, C.~Kreuzer, and R.~Stevenson}, {\em Instance optimality of the adaptive maximum strategy}, Found. Comput. Math., 16 (2016), pp.~33--68. The basic ingredients of our anlaysis consists of discrete reliability, discrete efficiency, overlay and closure estimates of local mesh refinement. The discrete efficiency requires a facet-interior node property of simplicial mesh refinement. The novel part of the optimality analysis is an elementary integral-sum argument that converts energy decrease into a cumulative cardinality bound of marked elements, yielding rate optimality without a contraction argument. Without the interior node property, we also prove that AFEM based on maximum marking and a single prescribed bisection per marked simplex is nearly rate-optimal up to a logarithmic factor. The key argument bounds the sum of error indicators on marked elements over a block of iterations by dyadically decomposing changes in the
Galerkin solution.
\end{abstract}

\begin{keywords}
adaptive finite element method, maximum marking, optimal convergence rate, a posteriori error estimation, interior node property, nearly optimal convergence rate
\end{keywords}

\begin{MSCcodes}
65N15, 65N30, 65N50
\end{MSCcodes}

\section{Introduction}\label{sec:introduction}

Adaptive finite element methods  (AFEMs) use a posteriori error indicators to
concentrate mesh refinement where it is needed. Their analysis asks
whether the computed solutions converge and whether the generated meshes
achieve the accuracy available on the best meshes of comparable size.
For elliptic boundary value problems, these questions have led to
convergence and optimality results for a variety of adaptive algorithms
\cite{BDD2004,CFPP2014,CKNS2008,Stevenson2007}. The marking rule is
central to this analysis: it determines which local error contributions
are reduced and how much refinement is required to reduce them.

Two common marking rules are D\"orfler bulk marking and maximum marking.
On the $\ell$-th adaptive mesh $\mathcal{T}_\ell$, let $\eta_\ell(T)$ denote
the local error indicator of element $T$. Bulk marking selects a minimal set
$\mathcal{M}_\ell^{\rm Dorf}\subset \mathcal{T}_\ell$ of elements carrying
a prescribed fraction of the global error estimator:
\[
\sum_{T\in\mathcal{M}_\ell^{\rm Dorf}}\eta_\ell(T)^2\geq\theta_{\rm Dorf}\sum_{T\in\mathcal{T}_\ell}\eta_\ell(T)^2,\qquad\theta_{\rm Dorf}\in(0,1).
\]
Maximum marking instead selects every
element whose indicator is at least a fixed fraction of the largest local
indicator. The latter rule is simple to implement and requires no
selection of a nearly minimal bulk set. Maximum marking has its roots in the pioneering work of Babuška and collaborators~\cite{BR1978,BV1984} and has been widely used in engineering community. Bulk marking was introduced later by Dörfler to establish theoretical convergence of AFEM.

The convergence of bulk marking algorithms is well understood.
For the two-dimensional Poisson problem, D\"orfler~\cite{Dorfler1996}
proved that AFEM with his bulk marking contracts the
energy error by a uniform factor at each step.
Optimal algebraic rates for AFEMs were proved in~\cite{BDD2004,CKNS2008,Stevenson2007}, and the
corresponding analysis was formulated through axioms of
adaptivity in~\cite{CFPP2014}. 
Recently, Diening and Storn~\cite{DS2026} removed the
smallness restriction on the bulk parameter,
proving rate optimality of AFEM with conforming Lagrange elements and
D\"orfler marking for every $\theta_{\rm Dorf}\in(0,1)$.
These rate results concern bulk marking
and do not directly establish optimal rates for the maximum marking rule.

Plain convergence without rate under general marking conditions that include maximum
marking was established in~\cite{MSV2008,GP2022}. 
Diening, Kreuzer, and Stevenson~\cite{DKS2016} proved instance optimality
for a modified maximum strategy for the two-dimensional Poisson problem.
Their strategy compares accumulated edge indicators over refinement
closures and removes contributions already covered by previous marks.
Instance optimality requires each generated mesh to compete with every
admissible mesh of sufficiently smaller cardinality. It is stronger
than rate optimality, which preserves each algebraic rate available in
a best-approximation class. The result in~\cite{DKS2016} therefore gives
a stronger performance guarantee, but for a different marking rule.
The modified strategy was extended to arbitrary fixed polynomial
degree in~\cite{IP2021}.

Currently, either instance or rate optimality of maximum marking strategy beyond two-dimension remains open.
In the present work, we prove rate optimality for standard maximum marking. For a marking
parameter $0<\theta\le1$, the marking set is exactly
\[
 \mathcal{M}_\ell=\left\{T\in\mathcal{T}_\ell:
 \eta_\ell(T)\ge\theta\max_{K\in\mathcal{T}_\ell}\eta_\ell(K)\right\}.
\]
We first present the proof for linear elements in two dimensions. We use edge-interior nodes to
control the jump terms in the discrete efficiency estimate, as in~\cite{DKS2016}.
A comparison mesh from the approximation class, together with overlay
and localized discrete reliability, gives a lower bound for the largest
current indicator. Maximum marking and discrete efficiency then bound the
decrease of the augmented energy from below. Integrating this estimate bounds
the cumulative number of marked elements and gives the optimal algebraic rate.
This short proof extends directly to arbitrary dimension $d\ge2$ by imposing
facet-interior nodes needed for discrete efficiency.
Our main result, Theorem~\ref{thm:p1-optimal}, has two
notable features. First, it establishes rate optimality for every maximum
marking parameter $0<\theta\le1$.
Second, the proof establishes optimal rates directly, bypassing the usual
route from contraction to optimality in analyses of D\"orfler marking.

A natural question is whether the additional bisections needed to create
facet-interior nodes are essential for rate optimality. When only one
bisection is prescribed per marked simplex, a facet carrying a large normal
derivative jump may survive conformity closure, so the preceding discrete
efficiency argument no longer applies at a single iteration. Surprisingly,
this obstruction can be overcome by
controlling the sum of the squared indicators on marked elements over blocks
of iterations rather than requiring discrete efficiency at every iteration.
The key observation is that standard estimator reduction together with a dyadic decomposition of the Galerkin increments gives the block estimate and yields the nearly optimal rate in
Theorem~\ref{thm:p1-one-bisection}, with only a logarithmic loss. To our
knowledge, this is the first nearly
optimal rate in terms of mesh cardinality for standard elementwise maximum
marking with one prescribed bisection per marked simplex.

Section~\ref{sec:model} presents the model problem and proves rate optimality
when the refinement step bisects each marked triangle and then bisects both
children. 
Section~\ref{sec:one-bisection} proves the nearly optimal rate in two
dimensions when only one bisection per marked triangle is prescribed.
Section~\ref{sec:setting} extends
both rate-optimality analysis  to simplicial meshes in arbitrary dimension using
Maubach bisection and discusses higher-order elements.

\section{Linear element in two dimension}\label{sec:model}
\subsection{Model problem and basic estimates}
Let $\Omega\subset\mathbb R^2$ be a bounded polygonal Lipschitz domain,
let $\mathcal{T}_0$ be a conforming initial triangulation, and let
$f\in L^2(\Omega)$. The exact solution $u\in H^1_0(\Omega)$ satisfies
\begin{equation}\label{eq:weak}
 a(u,v):=\int_\Omega\nabla u\cdot\nabla v\,\dd x
 =(f,v),\qquad v\in H^1_0(\Omega).
\end{equation}
For a conforming mesh $\mathcal{T}$, let $V(\mathcal{T})$ be the continuous
piecewise affine subspace of $H^1_0(\Omega)$ and define
$u_{\mathcal{T}}\in V(\mathcal{T})$ by
\begin{equation}\label{eq:discrete}
 a(u_{\mathcal{T}},v)=(f,v),\qquad v\in V(\mathcal{T}).
\end{equation}
Let $\mathbb T$
be the forest of conforming refinements of $\mathcal{T}_0$ based on newest vertex bisection (NVB) with a fixed initialization. These meshes are uniformly shape regular. We write
$\widetilde{\mathcal{T}}\succeq\mathcal{T}$ provided $\widetilde{\mathcal{T}}\in\mathbb{T}$ is a conforming refinement of $\mathcal{T}\in\mathbb{T}$ by newest vertex bisection, and
$\mathcal{T}\vee\mathcal{T}_*$ for the minimal common refinement of $\mathcal{T}$ and $\mathcal{T}_*$. 

For each triangle $T$ and edge $F$, set $h_T=|T|^{1/2}$
and $h_F={\rm diam}(F)$. Since $\Delta u_{\mathcal{T}}|_T=0$, the residual indicator is
\begin{equation}\label{eq:model-estimator}
 \eta_{\mathcal{T}}(T)^2=h_T^2\|f\|_{L^2(T)}^2+
 \frac12\sum_{\substack{F\subset\partial T\\F\text{ interior}}}h_F\|\llbracket\partial_nu_\mathcal{T}\rrbracket\|_{L^2(F)}^2.
\end{equation}
Here $\llbracket\partial_nu_\mathcal{T}\rrbracket$ is the jump of the normal derivative of $u_{\mathcal{T}}$ across $F$.
For $\mathcal{U}\subseteq\mathcal{T}$, we set
\begin{align*}
  \eta_{\mathcal{T}}(\mathcal U)^2&=\sum_{T\in\mathcal U}\eta_{\mathcal{T}}(T)^2,\qquad
\eta_{\mathcal{T}}=\eta_{\mathcal{T}}(\mathcal{T}),\\
\osc_\mathcal{T}(f)^2&=\sum_{T\in\mathcal T} h_T^2\|f-f_T\|_{L^2(T)}^2,\qquad f_T=\frac{1}{|T|}\int_T f\dd x.
\end{align*}
Reliability and efficiency of $\eta_\mathcal{T}$ reads
\begin{subequations}\label{eq:upper_lower}
  \begin{align}
    &|u-u_{\mathcal{T}}|_{H^1(\Omega)}\lesssim \eta_\mathcal{T},\\
    &\eta_\mathcal{T}\lesssim|u-u_{\mathcal{T}}|_{H^1(\Omega)}+{\rm osc}_\mathcal{T}(f).
  \end{align}
\end{subequations}
The symbols $\lesssim$ and $\eqsim$ denote bounds and equivalence up to
multiplicative constants independent of the mesh levels and the load $f$.
Subscripts indicate dependence on parameters such as the marking parameter $\theta$.

In Algorithm~\ref{alg:afem}, $m$ is the number of bisection generations
prescribed for each marked triangle. The first generation bisects the
triangle, and each subsequent generation bisects all of its current descendants.

\begin{algorithm}[htbp]
\caption{\textsf{AFEM}($m$) with maximum marking}
\label{alg:afem}
\begin{algorithmic}[1]
\Require Initial mesh $\mathcal{T}_0$, load $f\in L^2(\Omega)$,
         marking parameter $0<\theta\le1$, and integer $m\ge1$.
\For{$\ell=0,1,2,\ldots$}
  \State \textbf{SOLVE:} compute $u_\ell\in V(\mathcal{T}_\ell)$ satisfying
         \eqref{eq:discrete} on $\mathcal{T}_\ell$.
  \State \textbf{ESTIMATE:} compute $\eta_\ell(T)$ for all $T\in\mathcal{T}_\ell$.
  \State \textbf{MARK:} select the marking set 
\begin{equation}\label{eq:maximum}
 \mathcal{M}_\ell=\left\{T\in\mathcal{T}_\ell:
 \eta_\ell(T)\ge\theta\max_{K\in\mathcal{T}_\ell}\eta_\ell(K)\right\}.
\end{equation}
  \State \textbf{REFINE:} apply $m$ generations of NVB to each $T\in\mathcal{M}_\ell$,
         then restore conformity to obtain $\mathcal{T}_{\ell+1}$.
\EndFor
\end{algorithmic}
\end{algorithm}

\begin{figure}[ht]
\centering
\begin{tikzpicture}[scale=0.95,line width=0.75pt]
\coordinate (A) at (0,0);\coordinate (B) at (3.1,0);\coordinate (C) at (1.0,2.3);
\draw (A)--(B)--(C)--cycle;
\node[below] at (1.55,-0.18){a marked triangle};
\draw[-{Latex[length=3mm]},thick] (3.45,1.10)--(4.85,1.10);
\coordinate (a) at (5.2,0);\coordinate (b) at (8.3,0);\coordinate (c) at (6.2,2.3);
\coordinate (m) at (6.75,0);\coordinate (n) at (5.7,1.15);\coordinate (p) at (7.25,1.15);
\draw (a)--(b)--(c)--cycle;
\draw (c)--(m); \draw (n)--(m)--(p);
\node[below,align=center] at (6.75,-0.18){two generations of bisection};
\fill (m) circle (1.7pt);\fill (n) circle (1.7pt);\fill (p) circle (1.7pt);
\end{tikzpicture}
\caption{Two NVB generations split every original edge and produce edge interior nodes.}
\label{fig:twogen}
\end{figure}
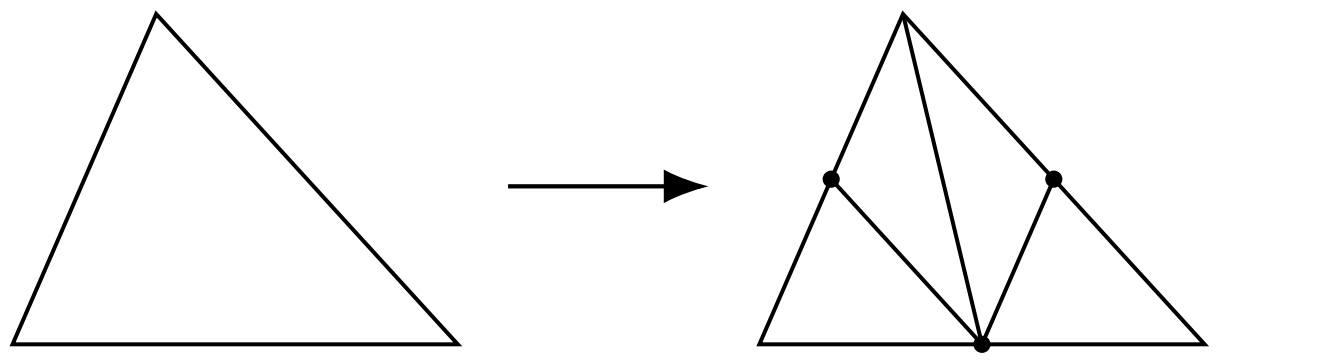

In this section we take $m=2$. Figure~\ref{fig:twogen} shows the
refinement of a marked triangle in \textsf{AFEM}($2$).

The two bisection generations create an interior node on each edge of every marked triangle.
Define the data term and augmented energy by
\begin{align*}
 D(\mathcal{U})=\sum_{T\in\mathcal{U}}h_T^2\|f\|_{L^2(T)}^2, \qquad E(\mathcal{T})=\en{u-u_{\mathcal{T}}}^2+D(\mathcal{T}).
\end{align*}
For $\widetilde{\mathcal{T}}\succeq\mathcal{T}$, every descendant of a refined triangle has at most half its area, so
\begin{equation}\label{eq:model-data-reduction}
 \tfrac12D(\mathcal{T}\setminus\widetilde{\mathcal{T}})\le D(\mathcal{T})-D(\widetilde{\mathcal{T}})\le D(\mathcal{T}\setminus\widetilde{\mathcal{T}}).
\end{equation}
Galerkin orthogonality therefore gives
\begin{equation}\label{eq:model-identity}
  \begin{aligned}
     E(\mathcal{T})-E(\widetilde{\mathcal{T}})
 &=|u_\mathcal{T}-u_{\widetilde{\mathcal{T}}}|^2_{H^1(\Omega)}+D(\mathcal{T})-D(\widetilde{\mathcal{T}})\\
 &\eqsim|u_\mathcal{T}-u_{\widetilde{\mathcal{T}}}|^2_{H^1(\Omega)}+D(\mathcal{T}\setminus\widetilde{\mathcal{T}})\ge0.
  \end{aligned}
\end{equation}

At level $\ell$, write $u_\ell=u_{\mathcal{T}_\ell}$, $\eta_\ell=\eta_{\mathcal{T}_\ell}$, and $E_\ell=E(\mathcal{T}_\ell)$.
For each interior edge $F$ of $\mathcal{T}_\ell$, let $\omega_\ell(F)$ be the union of the two triangles sharing $F$.
If $F$ is an edge of a marked triangle, its new interior node allows us to construct a discrete edge bubble in $V(\mathcal{T}_{\ell+1})$ and prove
\begin{equation}\label{eq:model-edge-test}
 h_F\|\llbracket\partial_nu_\ell\rrbracket\|_{L^2(F)}^2\lesssim|u_{\ell+1}-u_\ell|_{H^1(\omega_\ell(F))}^2
 +\sum_{K\subset\omega_\ell(F),\, K\in\mathcal{T}_\ell}|K|\|f\|_{L^2(K)}^2,
\end{equation}
see~\cite{MNS2000}.
Summing~\eqref{eq:model-edge-test} and using \eqref{eq:model-identity}, we obtain the discrete efficiency
\begin{equation}\label{eq:dis_lower}
 E_\ell-E_{\ell+1}\gtrsim\eta_\ell(\mathcal{M}_\ell)^2.
\end{equation}
The discrete reliability bound also holds~\cite{CKNS2008}
\begin{equation}\label{eq:dis_reliability}
  |u_{\widetilde{\mathcal{T}}}-u_{\mathcal{T}}|_{H^1(\Omega)}\lesssim \eta_\mathcal{T}(\mathcal{R}^+(\mathcal{T},\widetilde{\mathcal{T}}\,)),\qquad \forall\widetilde{\mathcal{T}}\succeq\mathcal{T}, 
\end{equation}
where $\mathcal{R}^+(\mathcal{T},\widetilde{\mathcal{T}})$ is the one-ring
enlargement of the set of refined triangles in $\mathcal{T}$:
\[
 \mathcal{R}^+(\mathcal{T},\widetilde{\mathcal{T}})=(\mathcal{T}\setminus\widetilde{\mathcal{T}})\cup
 \big\{K\in\mathcal{T}:K\text{ shares a face with some }T\in\mathcal{T}\setminus\widetilde{\mathcal{T}}\big\}.
\]
In addition, we need the following overlay and cumulative closure estimates of NVB (cf.~\cite{BDD2004,Stevenson2008,CKNS2008})  
\begin{subequations}
  \begin{align}
  \card(\mathcal{T}_1\vee\mathcal{T}_2)&\leq\card\mathcal{T}_1+\card\mathcal{T}_2-\card\mathcal{T}_0,\qquad\forall\mathcal{T}_1, \mathcal{T}_2\in\mathbb{T},\label{eq:overlay}\\
\card\mathcal{T}_L-\card\mathcal{T}_0
 &\lesssim\sum_{\ell=0}^{L-1}\card\mathcal{M}_\ell,\qquad\forall L\in\mathbb{N}_+.\label{eq:closure}
  \end{align}
\end{subequations}

\subsection{The optimality argument}\label{subsec:model-optimality}
For $s>0$, define the approximation seminorm
\begin{equation}\label{eq:approximation-class}
 |(u,f)|_{\A_s}:=\sup_{N\in\mathbb N_0}(N+1)^s
 \min_{\substack{\mathcal{T}\in\mathbb T\\\card\mathcal{T}-\card\mathcal{T}_0\le N}}
 E(\mathcal{T})^{1/2},
\end{equation}
which was introduced in \cite{BDD2004,CKNS2008}.
The following theorem establishes the rate optimality of \textsf{AFEM}($2$).

\begin{theorem}\label{thm:model-optimal}
  Let $s>0$ such that $|(u,f)|_{\A_s}<\infty$. Then for every $L\geq1$ and $0<\theta\le1$, \emph{\textsf{AFEM}(2)} satisfies
\begin{equation*}
|u-u_L|_{H^1(\Omega)}+D(\mathcal{T}_L)^{1/2}
 \le C_{s,\theta}|(u,f)|_{\A_s}
 (\card\mathcal{T}_L-\card\mathcal{T}_0)^{-s}.
\end{equation*}
\end{theorem}
\begin{proof}
Using the discrete efficiency \eqref{eq:dis_lower} and maximum marking criterion \eqref{eq:maximum}, we obtain a lower bound for energy decrease 
\begin{equation}\label{eq:marked-descent}
 E_\ell-E_{\ell+1}\gtrsim\theta^2\card\mathcal{M}_\ell
 \max_{T\in\mathcal{T}_\ell}\eta_\ell(T)^2.
\end{equation}
On the other hand, it follows from the discrete reliability \eqref{eq:dis_reliability} and the elementary cardinality estimate $\#\mathcal{R}^+(\mathcal{T},\widetilde{\mathcal{T}}\,)\lesssim\card\widetilde{\mathcal{T}}-\card\mathcal{T}$ that 
\begin{equation}\label{eq:comparison}
0\leq E(\mathcal{T})-E(\widetilde{\mathcal{T}})
 \lesssim(\card\widetilde{\mathcal{T}}-\card\mathcal{T})\max_{T\in\mathcal{T}}\eta_\mathcal{T}(T)^2.
\end{equation}

Write $A=|(u,f)|_{\A_s}$. Fix $L\ge1$ with $E_L>0$.
For each $\ell<L$, choose $N=\lceil(2A^2/E_\ell)^{1/(2s)}\rceil$.
A minimizing mesh $\mathcal{T}_*$ in~\eqref{eq:approximation-class} satisfies
\begin{equation}\label{eq:comp-card}
 E(\mathcal{T}_*)\le A^2(N+1)^{-2s}\le\tfrac12E_\ell,
 \qquad \card\mathcal{T}_*-\card\mathcal{T}_0\lesssim_s A^{1/s}E_\ell^{-1/(2s)}.
\end{equation}
Set $\widetilde{\mathcal{T}}=\mathcal{T}_*\vee\mathcal{T}_\ell$.
Using \eqref{eq:comp-card},
\eqref{eq:comparison},  the overlay estimate \eqref{eq:overlay}, we have 
\begin{align*}
   \tfrac12E_\ell&\le E_\ell-E(\mathcal{T}_*)\le E_\ell-E(\widetilde{\mathcal{T}})
 \\
 &\lesssim(\card\widetilde{\mathcal{T}}-\card\mathcal{T}_\ell)
 \max_{T\in\mathcal{T}_\ell}\eta_\ell(T)^2\\
 &\lesssim_s A^{1/s}E_\ell^{-1/(2s)}
 \max_{T\in\mathcal{T}_\ell}\eta_\ell(T)^2.
\end{align*}
Combining this lower bound for the maximum squared indicator with
\eqref{eq:marked-descent} gives
\begin{equation}\label{eq:descent}
 E_\ell-E_{\ell+1}
 \gtrsim_s\theta^2 A^{-1/s}\card\mathcal{M}_\ell E_\ell^{1+1/(2s)}.
\end{equation}
Since $t^{-1-1/(2s)}$ is decreasing, summing \eqref{eq:descent} yields
\begin{align*}
 \sum_{\ell=0}^{L-1}\card\mathcal{M}_\ell
 &\lesssim_{s,\theta}A^{1/s}
 \sum_{\ell=0}^{L-1}(E_\ell-E_{\ell+1})E_\ell^{-1-1/(2s)}\\
 &\le C_{s,\theta}A^{1/s}\int_{E_L}^{E_0}t^{-1-1/(2s)}\,\dd t\\
 &\le C_{s,\theta}A^{1/s}\big(E_L^{-1/(2s)}-E_0^{-1/(2s)}\big)\lesssim_{s,\theta}A^{1/s}E_L^{-1/(2s)}.
\end{align*}
Then by the cumulative closure estimate \eqref{eq:closure} and $A\ge E_0^{1/2}\ge E_L^{1/2}$, we arrive at  
\[
 \card\mathcal{T}_L-\card\mathcal{T}_0
 \lesssim_{s,\theta}A^{1/s}E_L^{-1/(2s)}.
\]
The proof is complete. 
\end{proof}

For longest-edge bisection in two dimensions,
the closure estimate~\eqref{eq:closure} was established
in~\cite{LiYangLEB2026}. Extending the present optimality proof
to this rule requires a separate argument because the overlay estimate~\eqref{eq:overlay} may not be true.

\section{Nearly optimal rates with one NVB bisection}\label{sec:one-bisection}
In this section, 
we return to  Algorithm~\ref{alg:afem} with $m=1$, that is, \textsf{AFEM}($1$) bisects each
marked triangle once along its reference edge before conformity closure.
An interior edge of a marked triangle can remain unsplit after closure,
so the per-iteration discrete efficiency estimate~\eqref{eq:dis_lower}
can fail. 

Remarkably, \textsf{AFEM}(1) is still nearly rate-optimal up in terms of mesh cardinality to a logarithmic loss. The key ingredient is the following  bound for the sum of marked
indicators over a block of consecutive adaptive iterations.
\begin{lemma}\label{lem:one-bisection-block}
For any integers $0\le L_1<L_2$,
\begin{equation*}
 \sum_{\ell=L_1}^{L_2-1}\eta_\ell(\mathcal M_\ell)^2
 \lesssim E_{L_2}+
 [1+\log(L_2-L_1)]^2(E_{L_1}-E_{L_2}).
\end{equation*}
\end{lemma}

\begin{proof}
For $v\in V(\mathcal T_\ell)$ and $g\in L^2(\Omega)$, we introduce
\[
 \widehat\eta_\ell(T;v,g)^2
 =h_T^2\|g\|_{L^2(T)}^2
 +\frac{h_T}{2}\sum_{\substack{F\subset\partial T\\F\text{ interior}}}
       \|\llbracket\partial_n v\rrbracket\|_{L^2(F)}^2
\]
to emphasize the estimator dependence on the finite element function $v$ and source $g$.
Set $\widehat\eta_\ell(\mathcal U;v,g)^2
=\sum_{T\in\mathcal U}\widehat\eta_\ell(T;v,g)^2$.
The element size $h_T$, rather than the facet size $h_F$ in the standard
estimator $\eta_\ell$, is essential below.
Since $h_T\eqsim h_F$ by shape regularity, \eqref{eq:model-estimator} gives
$\widehat\eta_\ell(\mathcal U;u_\ell,f)\eqsim\eta_\ell(\mathcal U)$.

Fix $L_1\le q<L_2$, $v\in V(\mathcal T_q)$, and $g\in L^2(\Omega)$.
For $\ell\ge q$, the jump of $\partial_n v$ vanishes on new edges inside
each $T\in\mathcal T_\ell$. On old edges, jump contributions are preserved
under subdivision, and $h_t\le 2^{-1/2}h_T$ for every descendant $t$ of a
bisected triangle $T$. Moreover, every marked triangle is bisected, so
$\mathcal M_\ell\subseteq\mathcal R_\ell:=\mathcal T_\ell\setminus\mathcal T_{\ell+1}$.
The standard estimator-reduction argument (cf.~\cite[Corollary~3.4]{CKNS2008}) therefore gives the following
inequalities:
\begin{align*}
 &(1-2^{-1/2})\widehat\eta_\ell(\mathcal M_\ell;v,g)^2\\
 &\quad\le (1-2^{-1/2})\widehat\eta_\ell(\mathcal R_\ell;v,g)^2\\
 &\quad\le \sum_{T\in\mathcal R_\ell}
   \left(\widehat\eta_\ell(T;v,g)^2
     -\sum_{\substack{t\in\mathcal T_{\ell+1}\\t\subsetneq T}}
        \widehat\eta_{\ell+1}(t;v,g)^2\right)\\
 &\quad= \widehat\eta_\ell(\mathcal T_\ell;v,g)^2
      -\widehat\eta_{\ell+1}(\mathcal T_{\ell+1};v,g)^2.
\end{align*}
Summing the preceding inequality over $\ell$ gives
\begin{equation}\label{eq:one-bisection-fixed}
 \sum_{\ell=q}^{L_2-1}\widehat\eta_\ell(\mathcal M_\ell;v,g)^2
 \lesssim \widehat\eta_q(\mathcal T_q;v,g)^2 .
\end{equation}

Set $K=1+\lfloor\log_2(L_2-L_1)\rfloor$. For $L_1\le\ell<L_2$, let
$\mathcal I(L_1,\ell)$ be the set of intervals that are
maximal under inclusion among the dyadic subintervals of
$\lbrack L_1,\ell\rparen$ of the form
\[
 \lbrack L_1+m2^r,\,L_1+(m+1)2^r\rparen,\qquad m,r\in\mathbb N_0.
\]
Members of $\mathcal I(L_1,\ell)$ are pairwise disjoint and form a  partition of $\lbrack L_1,\ell\rparen$. There are at most $K$ possible different lengths for the maximal intervals in $\mathcal I(L_1,\ell)$. At each fixed length $2^r$,  $\mathcal I(L_1,\ell)$ contains at most one
interval of this length, and $\card\mathcal I(L_1,\ell)\le K$.
Figure~\ref{fig:dyadic-partition}
illustrates this partition.

\begin{figure}[htbp]
\centering
\begin{tikzpicture}[x=0.60cm,y=0.68cm,font=\footnotesize,
  grid interval/.style={draw=black!45,line width=0.4pt},
  maximal interval/.style={draw=black!85,fill=blue!22,line width=0.8pt}]
  \foreach \endpoint in {0,15}
    \draw[densely dashed,black!55] (\endpoint,0) -- (\endpoint,4.95);
  \foreach \length/\row in {8/4.6,4/3.65,2/2.7,1/1.75}{
    \node[left] at (0,\row) {length $\length$};
    \foreach \start in {0,\length,...,15}
      \draw[grid interval] (\start,\row-0.24) rectangle ++(\length,0.48);
  }
  \foreach \start/\length/\row in {0/8/4.6,8/4/3.65,12/2/2.7,14/1/1.75}{
    \draw[maximal interval] (\start,\row-0.24) rectangle ++(\length,0.48);
    \draw[maximal interval] (\start,0) rectangle ++(\length,0.48);
  }
  \node[left] at (0,0.24) {$\mathcal I(0,15)$};
  \foreach \endpoint in {0,8,12,14,15,16}
    \node[below] at (\endpoint,-0.05) {$\endpoint$};
\end{tikzpicture}
\caption{Maximal dyadic intervals in $\lbrack L_1,\ell\rparen=\lbrack 0,15\rparen$ for $L_2=16$.
The shaded intervals have lengths $8,4,2,1$ from left to right.}
\label{fig:dyadic-partition}
\end{figure}
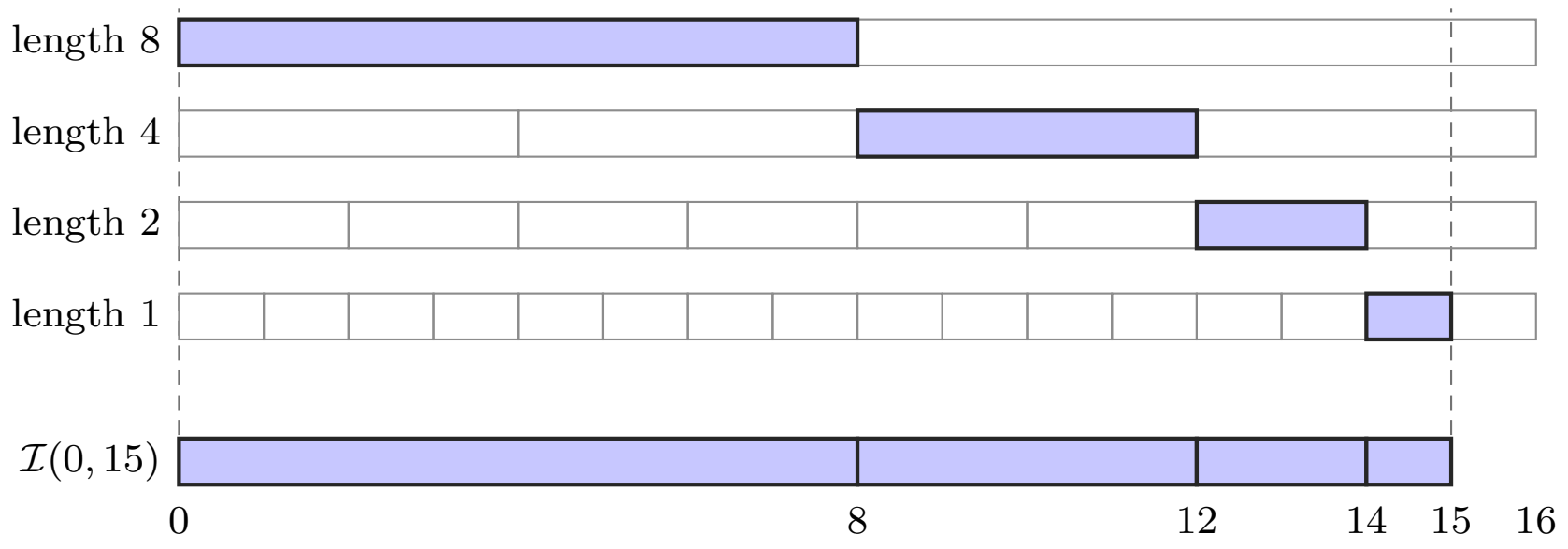

We make use of the telescope expansion $$u_\ell=u_{L_1}+\sum_{\lbrack p,q\rparen\in\mathcal I(L_1,\ell)}(u_q-u_p).$$
It follows from a triangle inequality for the seminorm
$(v,g)\mapsto\widehat\eta_\ell(\mathcal M_\ell;v,g)$ and Cauchy--Schwarz inequality that 
\begin{align*}
 &\widehat\eta_\ell(\mathcal M_\ell;u_\ell,f)^2\\
 &\quad\le 2\widehat\eta_\ell(\mathcal M_\ell;u_{L_1},f)^2
       +2\left(\sum_{\lbrack p,q\rparen\in\mathcal I(L_1,\ell)}
         \widehat\eta_\ell(\mathcal M_\ell;u_q-u_p,0)\right)^2\\
 &\quad\le 2\widehat\eta_\ell(\mathcal M_\ell;u_{L_1},f)^2
       +2\card\mathcal I(L_1,\ell)\sum_{\lbrack p,q\rparen\in\mathcal I(L_1,\ell)}
         \widehat\eta_\ell(\mathcal M_\ell;u_q-u_p,0)^2\\
 &\quad\le 2\widehat\eta_\ell(\mathcal M_\ell;u_{L_1},f)^2
       +2K\sum_{\lbrack p,q\rparen\in\mathcal I(L_1,\ell)}
         \widehat\eta_\ell(\mathcal M_\ell;u_q-u_p,0)^2 .
\end{align*}
Let $\mathcal I_r$ be the set of distinct intervals of length $2^r$ in
these partitions. Each $\lbrack p,q\rparen$ contributes only at levels $\ell\ge q$.
Applying \eqref{eq:one-bisection-fixed} first to $(u_{L_1},f)$ and then to
each $(u_q-u_p,0)\in V(\mathcal T_q)\times L^2(\Omega)$ gives
\begin{align*}
 \sum_{\ell=L_1}^{L_2-1}\eta_\ell(\mathcal M_\ell)^2
 &\lesssim \eta_{L_1}^2+K\sum_{r=0}^{K-1}\sum_{\lbrack p,q\rparen\in\mathcal I_r}
       \widehat\eta_q(\mathcal T_q;u_q-u_p,0)^2\\
 &\lesssim \eta_{L_1}^2+K\sum_{r=0}^{K-1}\sum_{\lbrack p,q\rparen\in\mathcal I_r}
       \en{u_q-u_p}^2 ,
\end{align*}
where the last inequality is the inverse trace estimate on $\mathcal T_q$.
At each length the distinct dyadic intervals are disjoint, so Galerkin
orthogonality~\eqref{eq:model-identity} yields
\[
 \sum_{\lbrack p,q\rparen\in\mathcal I_r}\en{u_q-u_p}^2
 \le \sum_{\lbrack p,q\rparen\in\mathcal I_r}(E_p-E_q)\le E_{L_1}-E_{L_2}.
\]
Finally, estimator efficiency~\eqref{eq:upper_lower} and
$\osc_{\mathcal T_{L_1}}(f)^2\le D(\mathcal T_{L_1})$ give
$\eta_{L_1}^2\lesssim E_{L_1}$. Hence
\[
 \sum_{\ell=L_1}^{L_2-1}\eta_\ell(\mathcal M_\ell)^2
 \lesssim E_{L_1}+K^2(E_{L_1}-E_{L_2})
 \lesssim E_{L_2}+K^2(E_{L_1}-E_{L_2}).
\]
The assertion follows from $K\lesssim 1+\log(L_2-L_1)$.
\end{proof}

The next theorem shows that the AFEM with standard maximum marking and one bisection for marked elements is rate-optimal up to a logarithmic factor.
\begin{theorem}\label{thm:one-bisection-rate}
Let $s>0$ and $|(u,f)|_{\A_s}<\infty$. For $0<\theta\le1$ and every integer $L\ge1$, \emph{\textsf{AFEM}(1)} satisfies
\begin{equation*}
 |u-u_L|_{H^1(\Omega)}+D(\mathcal T_L)^{1/2}
 \le C_{s,\theta}|(u,f)|_{\A_s}\,(\card\mathcal T_L-
\card\mathcal T_0)^{-s}(1+\log L)^{2s}.
\end{equation*}
\end{theorem}

\begin{proof}
If $E_L=0$, the statement is trivial.
Assume $E_L>0$ and set $A=|(u,f)|_{\A_s}$. Then $A\ge E_0^{1/2}\ge E_L^{1/2}>0$.
The comparison argument in Theorem \ref{thm:model-optimal} using
\eqref{eq:comparison}, \eqref{eq:comp-card}, and \eqref{eq:overlay}
does not rely on the two-generation REFINE step in \textsf{AFEM}(2). Therefore, we still have 
$$\max_{T\in\mathcal T_\ell}\eta_\ell(T)^2\gtrsim_s A^{-1/s}E_\ell^{1+1/(2s)}\quad\text{ for } \ell<L.$$ By the maximum marking 
\eqref{eq:maximum},
\begin{equation}\label{eq:one-bisection-marked-lower}
 \eta_\ell(\mathcal M_\ell)^2
 \gtrsim_s \theta^2 A^{-1/s}\,\card\mathcal M_\ell\, E_\ell^{1+1/(2s)} .
\end{equation}
For $k\ge0$ let $g_k=2^{-k}E_0$ and group the indices $\ell<L$ for
which $g_k/2<E_\ell\le g_k$. By monotonicity each nonempty group
is an interval of integers, denoted by $\lbrack L_1,L_2\rparen$ for this $k$.
Summing \eqref{eq:one-bisection-marked-lower} after solving for $\card\mathcal M_\ell$,
using $E_\ell>g_k/2$ on this group, and applying Lemma~\ref{lem:one-bisection-block} give
\begin{align*}
 \sum_{\ell=L_1}^{L_2-1}\card\mathcal M_\ell
 &\lesssim_{s,\theta} A^{1/s}\sum_{\ell=L_1}^{L_2-1}
   \frac{\eta_\ell(\mathcal M_\ell)^2}{E_\ell^{1+1/(2s)}}\\
 &\lesssim_s A^{1/s}g_k^{-1-1/(2s)}
   \sum_{\ell=L_1}^{L_2-1}\eta_\ell(\mathcal M_\ell)^2\\
 &\lesssim A^{1/s}g_k^{-1-1/(2s)}
   \bigl[E_{L_2}+(1+\log(L_2-L_1))^2(E_{L_1}-E_{L_2})\bigr]\\
 &\le A^{1/s}g_k^{-1/(2s)}(1+\log L)^2,
\end{align*}
where the last line uses $E_{L_2}\le E_{L_1}\le g_k$ and $L_2-L_1\le L$.
Let $K_L=\max\{k\in\mathbb N_0:g_k\ge E_L\}$. The cumulative closure
estimate~\eqref{eq:closure} and the preceding bound for each nonempty group
give
\begin{align*}
 \card\mathcal T_L-\card\mathcal T_0
 &\lesssim \sum_{\ell=0}^{L-1}\card\mathcal M_\ell\\
 &\lesssim_{s,\theta} A^{1/s}(1+\log L)^2
   \sum_{k=0}^{K_L}g_k^{-1/(2s)}\\
 &\lesssim_{s,\theta} A^{1/s}E_L^{-1/(2s)}(1+\log L)^2.
\end{align*}
Indeed, $g_k=2^{-k}E_0$, so the geometric sum is dominated by its last
term, and $g_{K_L}\ge E_L$. The proof is complete. 
\end{proof}

From Theorem \ref{thm:one-bisection-rate} with $L\leq \card\mathcal T_L-
\card\mathcal T_0$, we know the convergence rate of \textsf{AFEM}(1) is $O(N^{-s}\log(1+N)^{2s})$ in terms of the number $N$ of degrees of freedom.

\section{Optimality in arbitrary dimension}\label{sec:setting}
We extend the linear-element analysis of Section~\ref{sec:model} to
any fixed dimension $d\ge2$. Let $\Omega\subset\mathbb R^d$ be a
bounded polyhedral Lipschitz domain and retain $f\in L^2(\Omega)$.
The weak problem, Galerkin equations, and notation for the adaptive
sequence remain those of Section~\ref{sec:model}. In this section,
$V(\mathcal{T})$ is the continuous piecewise affine subspace of
$H^1_0(\Omega)$ on a simplicial mesh $\mathcal{T}$.

\subsection{Mesh refinement}
Let $\mathcal{T}_0$ be a conforming simplicial triangulation of $\Omega$. Let $\mathbb{T}$ be the forest of conforming simplicial refinement of $\mathcal{T}_0$ based on tagged Maubach bisection rule~\cite{Maubach1995}, with a
compatible initialization on special $\mathcal{T}_0$ \cite{Stevenson2008} or the Alkamper--Gaspoz--Klofkorn, Diening--Gehring--Storn initialization that works on arbitrary initial mesh $\mathcal{T}_0$ \cite{AGK2018,DGS2026}.
The notation $\widetilde{\mathcal{T}}\succeq\mathcal{T}$
retains its meaning from Section~\ref{sec:model}. The error indicator is still~\eqref{eq:model-estimator}, with $h_T=|T|^{1/d}$ and $F$ now a
$(d-1)$-dimensional facet. We retain the estimator and adaptive-level
notation of Section~\ref{sec:model}.

The proof of Theorem~\ref{thm:model-optimal} is based on four ingredients: the discrete efficiency \eqref{eq:dis_lower}, discrete reliability \eqref{eq:dis_reliability}, the overlay estimate \eqref{eq:overlay}, and the cumulative closure estimate \eqref{eq:closure}. Clearly, \eqref{eq:dis_reliability}--\eqref{eq:closure} are still available in arbitrary dimension. Regarding the discrete efficiency, it suffices
to bisect marked simplices in $\mathcal{M}_\ell$ sufficiently many times so that every facet of $\mathcal{M}_\ell$ contains an interior node, and thus \eqref{eq:dis_lower} holds true.

For $\mathcal U\subseteq\mathcal{T}$, let
$\operatorname{Refine}_q(\mathcal{T},\mathcal U)$ be the coarsest conforming
refinement obtained by $q$ generations of Maubach bisection of every
$T\in\mathcal U$ and conformity closure.
A generation means one bisection of every current descendant of an
original element. As for the two-generation NVB step in Section~\ref{sec:model},
fixed $q$ only changes the constant in the standard cumulative closure
estimate~\eqref{eq:closure}.

\begin{lemma}\label{lem:resolution}
For a simplex $T$ in any dimension $d\ge2$, there is a fixed
integer $q>0$ such that $q$ generations of Maubach bisection create a vertex in the relative
interior of every facet of $T$.
\end{lemma}

\begin{figure}[htbp]
\centering
\includegraphics[width=\linewidth]{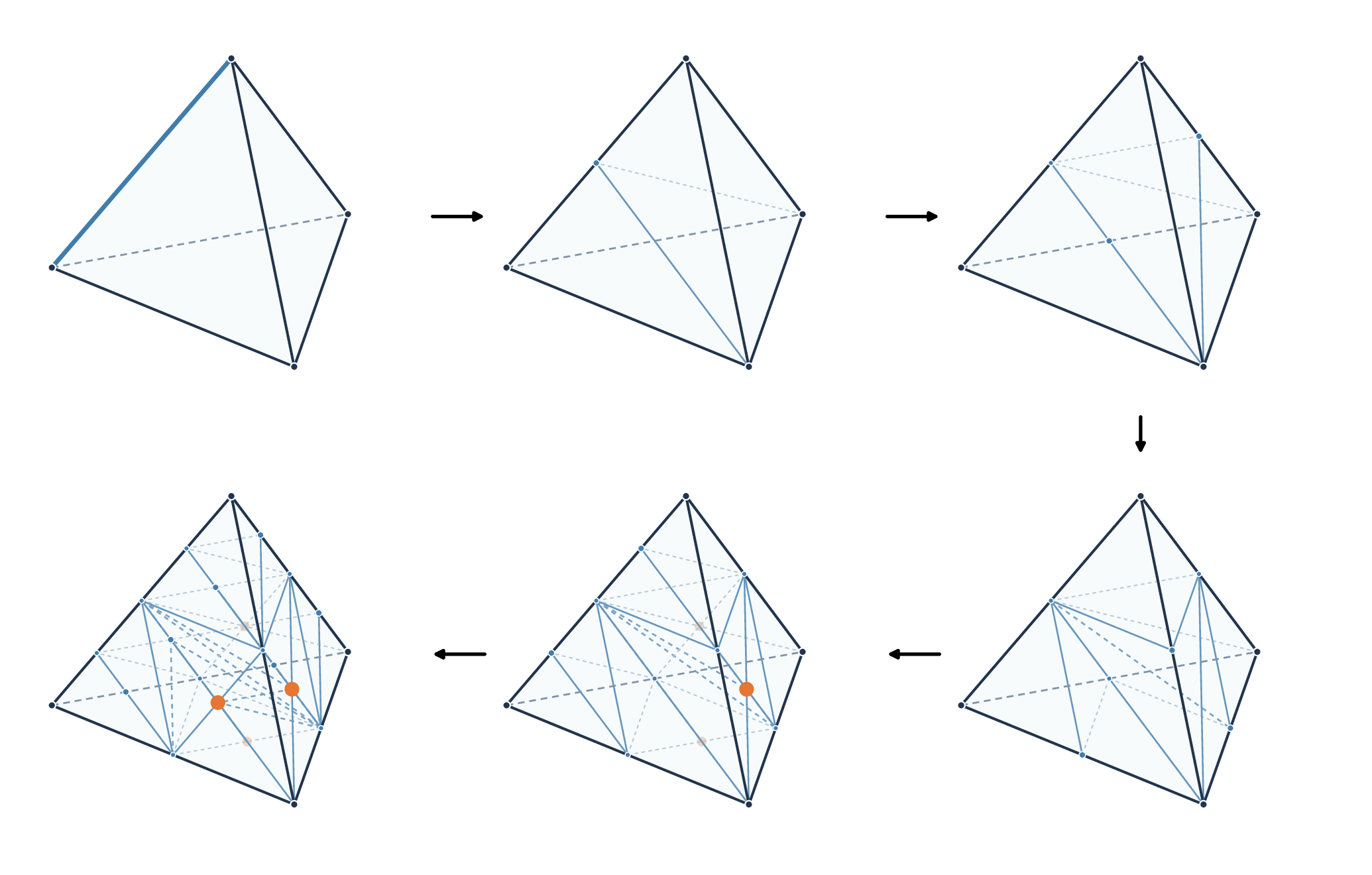}
\caption{Five successive Maubach bisection generations of a tetrahedron
 with initial tag $3$. The final refinement contains $32$
 tetrahedra before conformity closure.}
\label{fig:tetra-resolution}
\end{figure}

Let $q$ be the minimal integer stated in Lemma \ref{lem:resolution}. In particular, $q=2$ for $d=2$ and $q=5$
for $d=3$. Figure~\ref{fig:tetra-resolution} shows the boundary
of an actual five-generation bisection in 3D. In this section we use
\textsf{AFEM}($q$) with Maubach bisection, so the refinement step gives
$\mathcal{T}_{\ell+1}=\operatorname{Refine}_q(\mathcal{T}_\ell,\mathcal{M}_\ell).$
Every original interior facet $F$ of a marked element then contains a
vertex $z_F\in F$ satisfying
\begin{equation}\label{eq:facet-resolution}
 \operatorname{dist}_F(z_F,\partial F)\gtrsim h_F.
\end{equation}
The SOLVE, ESTIMATE, and MARK steps are unchanged.

For the Maubach families of simplicial mesh refinement algorithms, the cumulative closure estimate have been established in \cite{Stevenson2008} for special $\mathcal{T}_0$, and recently even on arbitrary initial mesh \cite{Gehring2023,DGS2026,LiYangAMP2026}, respectively.
Overlay is given by the union of refinement forests, see~\cite[Corollary 4.18]{Gehring2023}.

\subsection{Discrete efficiency}
We retain $D$ and $E$ from Section~\ref{sec:model}, with
$h_T=|T|^{1/d}$. 
For $\widetilde{\mathcal{T}}\succeq\mathcal{T}$, every descendant of a refined element has at most half its volume,
so $h_t^2\le2^{-2/d}h_T^2$ for $t\subsetneq T$. Consequently,
\begin{equation*}
 (1-2^{-2/d})D(\mathcal{T}\setminus\widetilde{\mathcal{T}})\le D(\mathcal{T})-D(\widetilde{\mathcal{T}})
 \le D(\mathcal{T}\setminus\widetilde{\mathcal{T}}).
\end{equation*}
Galerkin orthogonality gives the same identity as in
\eqref{eq:model-identity}, now with
\begin{equation}\label{eq:p1-identity}
 E(\mathcal{T})-E(\widetilde{\mathcal{T}})
  \eqsim|u_\mathcal{T}-u_{\widetilde{\mathcal{T}}}|_{H^1(\Omega)}^2+D(\mathcal{T}\setminus\widetilde{\mathcal{T}})\ge0.
\end{equation}

We next prove discrete efficiency for this refinement step.
For an interior facet $F$ of a marked element, let $\phi_F\in V(\mathcal{T}_{\ell+1})$ be the hat function at $z_F$ on $\mathcal{T}_{\ell+1}$, which is
supported in $\omega_F$, the union of the two simplices in $\mathcal{T}_\ell$ sharing $F$.
Since the number of refinement generations, including conformity closure, is bounded in terms of the fixed \(q\), a scaling argument gives
\begin{equation}\label{eq:phiF}
 0\le\phi_F\le1,\qquad
 \|\nabla\phi_F\|_{L^\infty(\omega_F)}\lesssim h_F^{-1},\qquad\int_F\phi_F\,\dd s\gtrsim |F|.
\end{equation}

\begin{lemma}\label{lem:p1-eff}
It holds that 
\begin{equation*}
 \eta_\ell(\mathcal{M}_\ell)^2
 \lesssim E(\mathcal{T}_\ell)-E(\mathcal{T}_{\ell+1}).
\end{equation*}
\end{lemma}
\begin{proof}
Set $\delta=u_{\ell+1}-u_\ell$ and
$\mathcal{R}_\ell=\mathcal{T}_\ell\setminus\mathcal{T}_{\ell+1}$.
For each interior facet $F$ of a marked element, the jump
$J_F=\llbracket\partial_n u_\ell\rrbracket|_F$ is constant.
Use $v_F=J_F\phi_F\in V(\mathcal{T}_{\ell+1})$ as a test function.
The hat-function bounds give
\[
 \|v_F\|_{L^2(\omega_F)}+h_F\|\nabla v_F\|_{L^2(\omega_F)}
 \lesssim h_F^{1/2}\|J_F\|_{L^2(F)}.
\]
By~\eqref{eq:phiF}, the Galerkin equation and
integration by parts,
\[
 \|J_F\|_{L^2(F)}^2
 \lesssim\int_F\phi_FJ_F^2\,\dd s
 =(f,v_F)_{\omega_F}-(\nabla\delta,\nabla v_F)_{\omega_F}.
\]
Cauchy--Schwarz and $h_T\eqsim h_F$ on the patch imply
\begin{equation}\label{eq:p1-facet-test}
 h_F\|J_F\|_{L^2(F)}^2
 \lesssim\|\nabla\delta\|_{L^2(\omega_F)}^2
       +\sum_{T\subset\omega_F,\,T\in\mathcal{T}_\ell}h_T^2\|f\|_{L^2(T)}^2.
\end{equation}
Every marked element is refined. The new vertex inside $F$ also
forces both coarse neighbors of $F$ to be refined, so they belong
to $\mathcal{R}_\ell$. Summing~\eqref{eq:p1-facet-test}, using
bounded overlap of the facet patches, and adding the volume
contributions gives
\[
 \eta_\ell(\mathcal{M}_\ell)^2
 \lesssim\|\nabla\delta\|_{L^2(\Omega)}^2+D(\mathcal{R}_\ell).
\]
Combining this result with \eqref{eq:p1-identity} proves the assertion.
\end{proof}

We have verified discrete efficiency in arbitrary dimension. Therefore, w
The proof of Theorem~\ref{thm:model-optimal} now gives the following result.
\begin{theorem}\label{thm:p1-optimal}
If $|(u,f)|_{\A_s}<\infty$ for some $s>0$, then for every $0<\theta\le1$ and every integer  $L\ge1$, \emph{\textsf{AFEM}($q$)} satisfies
\begin{equation*}
|u-u_L|_{H^1(\Omega)}+D(\mathcal{T}_L)^{1/2}
 \le C_{s,\theta}|(u,f)|_{\A_s}
       (\card\mathcal{T}_L-\card\mathcal{T}_0)^{-s}.
\end{equation*}
\end{theorem}

Finally, we return to \textsf{AFEM}(1) with one prescribed Maubach bisection per
marked simplex. Extend $\widehat\eta_\ell(\mathcal U;v,g)$ from the proof of
Lemma~\ref{lem:one-bisection-block} to simplicial meshes in $\mathbb R^d$ by
letting $F$ range over the facets of $T$ and taking $h_T=|T|^{1/d}$.
Fix $v\in V(\mathcal T_j)$, $g\in L^2(\Omega)$, and $\ell\ge j$. Every
descendant $t$ of a bisected simplex $T$ satisfies
$h_t\le 2^{-1/d}h_T$. The same estimator-reduction argument gives
\[
 (1-2^{-1/d})\widehat\eta_\ell(\mathcal M_\ell;v,g)^2
 \le \widehat\eta_\ell(\mathcal T_\ell;v,g)^2
       -\widehat\eta_{\ell+1}(\mathcal T_{\ell+1};v,g)^2.
\]
Summing this inequality over $\ell$ gives the $d$-dimensional analogue
of~\eqref{eq:one-bisection-fixed}. The same dyadic summation, using the
inverse trace estimate on simplices and Galerkin
orthogonality~\eqref{eq:p1-identity}, gives the block estimate of
Lemma~\ref{lem:one-bisection-block} in dimension $d$. Apply the proof of Theorem~\ref{thm:one-bisection-rate} with the
$d$-dimensional block estimate establishes the next theorem.

\begin{theorem}\label{thm:p1-one-bisection}
Let $d\ge2$  and $|(u,f)|_{\A_s}<\infty$ for some $s>0$. Then for $0<\theta\le1$ and $L\ge1$,  \emph{\textsf{AFEM}(1)} satisfies 
\begin{equation*}
 |u-u_L|_{H^1(\Omega)}+D(\mathcal T_L)^{1/2}
 \le C_{s,\theta}|(u,f)|_{\A_s}\,(\card\mathcal T_L-
\card\mathcal T_0)^{-s}(1+\log L)^{2s}.
\end{equation*}
\end{theorem}

\begin{remark}[Higher-order elements]\label{rem:higher-order}
For any fixed $k\ge2$, use conforming $P_k$ elements with the standard
residual estimator. Let $\Pi_T^{k-2}$ be the $L^2(T)$-orthogonal
projection onto $\mathbb P_{k-2}(T)$ and set
\begin{align*}
 D(\mathcal{T})&=\sum_{T\in\mathcal{T}}h_T^2
       \|(1-\Pi_T^{k-2})f\|_{L^2(T)}^2,\\
 E(\mathcal{T})&=\en{u-u_{\mathcal{T}}}^2+D(\mathcal{T}).
\end{align*}
Let $q\in\mathbb{N}$ for bisection generation be large enough, such that $\operatorname{Refine}_q(\mathcal{T}_\ell,\mathcal{M}_\ell)$ produces facet-interior and element-interor nodes for all facets and elements in $\mathcal{M}_\ell$.
The resulting facet- and element-interior nodes allow the
projected-data estimates and higher-order residual arguments of ~\cite{IP2021}
to extend to any fixed dimension $d\ge2$, using discrete facet and element bubble tests. Combining these estimates with
the proof of Theorem~\ref{thm:model-optimal} yields the rate bound, with the present energy
in~\eqref{eq:approximation-class} and constants that may depend on $k$.

The nearly optimal bound in Theorem~\ref{thm:p1-one-bisection} remains valid
for \textsf{AFEM}(1) and every polynomial degree $k\ge2$. To see this, consider the proof of
Lemma~\ref{lem:one-bisection-block}. For a fixed pair $(v,g)$, the element
residual on $T$ is $g+\Delta v$. Under refinement, this residual is simply
restricted to the children of $T$, while $\partial_n v$ has no jump across
new facets inside $T$. Consequently, the estimator reduction
in~\cite[Corollary~3.4]{CKNS2008} yields~\eqref{eq:one-bisection-fixed}. For a
Galerkin increment $w=u_q-u_p$, standard polynomial inverse and trace
estimates give
$\widehat\eta_q(\mathcal T_q;w,0)\lesssim_k\|\nabla w\|_{L^2(\Omega)}$.
Moreover, the projected data term $D$ decreases under refinement by local
best approximation. The same dyadic summation then
yield the asserted rate with the logarithmic factor $(1+\log L)^{2s}$.
\end{remark}

\section*{AI Declarations}
During the preparation of this work, GPT-6.0 Astra (OpenAI)
was used to derive proof strategies of main theorems in this paper. The
author have reviewed, verified and rewritten the AI-assisted proofs. The author assumes full responsibility for
the results and the final manuscript.

\end{document}